\documentclass[12pt]{article}
\usepackage[cp1251]{inputenc}
\usepackage[english]{babel}
\usepackage{amssymb}
\usepackage{amsmath}
\usepackage{amsthm}
\usepackage{amsbsy}
\usepackage{amsfonts}
\fontsize{16}{20}
\newtheorem{thm}{Theorem}[section]

\newtheorem{prop}[thm]{Proposition}

\newtheorem{rem}[thm]{Remark}

\numberwithin{equation}{section}

\begin{document}
\begin{center}
\begin{Large}
{\bf A method for computing waveguide scattering matrices in the
presence of discrete spectrum}
\end{Large}

\vspace{4ex}

\begin{large}
B. A. Plamenevskii, O. V. Sarafanov
\end{large}
\footnote{The work of both authors was supported by Russian
Foundation for Basic Research project 09-01-00191-a}

\vspace{6ex}

\end{center}
\begin{abstract}
A waveguide $G$ lies in $\mathbb R^{n+1}$, $n\geq 1$, and outside
a large ball coincides with the union of finitely many
non-overlapping semi-cylinders ("cylindrical ends"). The waveguide
is described by the operator $\{\mathcal{L}(x, D_x)-\mu,
\mathcal{B} (x, D_x)\}$ of an elliptic boundary value problem in
$G$, where $\mathcal{L}$ is a matrix differential operator,
$\mathcal{B}$ is a boundary operator, and $\mu$ is a spectral
parameter. The operator $\{\mathcal{L}, \mathcal{B}\}$ is
self-adjoint with respect to a Green formula. The role of
$\mathcal{L}$ can be played, e.g., by the Helmholtz operator, by
the operators in elasticity theory and hydrodynamics. As
approximation for a row of the scattering matrix $S(\mu)$,  we
take the minimizer of a quadratic functional $J^R(\cdot, \mu)$. To
construct the functional, we solve an auxiliary boundary value
problem in the bounded domain obtained by truncating the
cylindrical ends of the waveguide at distance $R$. As $R\to
\infty$, the minimizer  $a (R, \mu)$ tends with exponential rate
to the corresponding row of the scattering matrix uniformly  on
every finite closed interval of the continuous spectrum not
containing the thresholds. Such an interval may contain
eigenvalues of the waveguide with eigenfunctions exponentially
decaying at infinity ("trapped modes"). Eigenvalues of this sort,
as a rule, occur in waveguides of complicated geometry. Therefore,
in applications, the possibility to avoid worrying about (probably
not detected) trapped modes turns out to be an important advantage
of the method.
\end{abstract}

\section{Introduction}\label{s1}

The waveguide to be considered in the paper occupies a domain $G$ in
$\mathbb{R}^{n+1}$ having several cylindrical outlets to infinity
("cylindrical ends"). This means that outside a large ball centered
at the origin the domain $G$ coincides with the union of
non-overlapping  semicylinders $\Pi_+^1, \dots,\Pi^P_+$; here
$\Pi^p_+=\{(y^p, t^p): y^p\in \Omega^p, t^p>0\}$, $(y^p, t^p)$ are
local coordinates in $\Pi^p_+$, and the cross-section $\Omega^p$ of
the cylinder $\Pi^p$ is a bounded $n$-dimensional domain with smooth
boundary $\partial \Omega^p$. The waveguide is described by the
operator $\{\mathcal{L}(x, D_x)-\mu, \mathcal{B} (x, D_x)\}$ of an
elliptic boundary value problem in  $G$, where $\mathcal{L}$ is a
matrix differential operator, $\mathcal{B}$ is a boundary condition
operator, and $\mu$ is a spectral parameter. The operator
$\{\mathcal{L}, \mathcal{B}\}$ is self-adjoint with respect to a
Green formula. At infinity, the coefficients of $\mathcal{L}$ and
$\mathcal{B}$ stabilize with exponential rate in every semicylinder
$\Pi^p_+$ to functions independent of the axial variable.

Let us consider the homogeneous problem
\begin{eqnarray}\label{1.1}
(\mathcal{L}(x, D_x)-\mu) u(x)=0, \qquad x\in G, \\
\nonumber \mathcal{B} (x, D_x)u(x)=0, \qquad x\in \partial G.
\end{eqnarray}
We assume that the interval $[\mu_1, \mu_2]\subset \mathbb{R}$
belongs to the continuous spectrum of the operator $\{\mathcal{L} -
\mu, \mathcal{B}\}$ and contains no threshold values of the spectral
parameter.  In other words, for every  $\mu \in [\mu_1, \mu_2]$
there exists the same (finite) number of solutions to the problem
(\ref{1.1}) linearly independent modulo $L_2 (G)$; such solutions
are called eigenfunctions of the continuous spectrum. The interval
$[\mu_1, \mu_2]$ may contain eigenvalues of the problem (\ref{1.1})
with eigenfunctions in  $L_2 (G)$. Any eigenfunction in
 $L_2 (G)$ is exponentially decaying at infinity while the eigenvalues
are of finite multiplicity and can not accumulate in   $[\mu_1,
\mu_2]$. Thus, when $\mu \in [\mu_1, \mu_2]$ turns out to be an
eigenvalue, the number of bounded solutions linearly independent in
the ordinary sense increases and yet the number of solutions
 linearly independent modulo $L_2 (G)$ (or, equivalently, modulo
exponentially decaying terms) remains constant on $[\mu_1, \mu_2]$;
denote the number by $M$. For any  $\mu \in [\mu_1, \mu_2]$ in the
space of continuous spectrum eigenfunctions there exists a basis
 $Y_1 (\cdot, \mu), \dots, Y_M (\cdot, \mu)$ modulo $L_2
(G)$ such that
$$
Y_j (x, \mu)=u^+_j (x, \mu)+\sum_{k=1}^{M}S_{jk}(\mu)u^-_k(x, \mu) +
O (e^{-\varepsilon |x|})
$$
for $|x|\to \infty$ and  $j=1, \dots, M$; here $\varepsilon$ is a
sufficiently small positive number,   $u^+_j (\cdot, \mu)$ are
incoming waves while $u^-_j (\cdot, \mu)$ are outgoing ones (precise
definitions see in 2.2). The matrix  $S(\mu)=\|S_{jk}(\mu)\|$ is
unitary; it is called the scattering matrix.

The paper is devoted to justification of an approximate computation
method for the scattering matrix. A detailed description of the
method has been given in \ref{ss2.3}. In brief as an approximation
to the $l$-th row $S_l (\mu)=(S_{l,1}(\mu), \dots, S_{l, M}(\mu))$
of the scattering matrix we take the minimizer  $a (R, \mu)$ of a
quadratic functional $J^R_l(\cdot, \mu)$. The functional is
constructed by solving an auxiliary boundary value problem in the
bounded domain $G^R$ obtained from $G$ by cutting off the
cylindrical ends at a sufficiently large distance $R$ from the
origin. In the present paper, we prove that for  $R\geqslant R_0$
and all $\mu \in [\mu_1, \mu_2]$ there holds the estimate
\begin{equation}\label{1.2}
 \|a (R, \mu)-S_l (\mu)\|\leqslant
Ce^{-\Lambda R}
\end{equation}
with some positive numbers $\Lambda$ and $C$ that are independent of
$R$ and $\mu$. Thus, as $R\to \infty$, the minimizer $a(R, \mu)$
tends to the corresponding row of the scattering matrix with
exponential rate uniformly with respect to $\mu \in [\mu_1, \mu_2]$.

As mentioned above, the interval $[\mu_1, \mu_2]$ of the continuous
spectrum may contain eigenvalues of the operator $\{\mathcal{L}
-\mu, \mathcal{B}\}$. In nonhomogeneous waveguides of complicated
geometry, as a rule, there occur trapped modes, that is,
eigenfunctions exponentially decaying at infinity. Therefore, in
applications,  the possibility to avoid worrying about (probably not
detected) trapped modes turns out to be an important advantage of
the method.

For the Helmholtz operator in a close situation, the method under
discussion was suggested in  \cite{GHNP}. The justification of the
method in \cite{GHNP} made use of Proposition 3 there (presented
without proof), which is valid only if the interval $[\mu_1, \mu_2]$
is free from the eigenvalues of the waveguide; this restriction was
not indicated in Proposition 3. (In our context the restriction
would mean that $[\mu_1, \mu_2]$ must be free from the eigenvalues
of the operator $\{\mathcal{L} -\mu, \mathcal{B}\}$ in $G$.)  In the
presence of waveguide eigenvalues in $[\mu_1, \mu_2]$, the
inequality (\ref{1.2}) was proved for the Helmholtz operator in
\cite{PS}. The proof given in the present paper for the general
elliptic problems is new for the Helmholtz operator as well and is
simpler than that in \cite{PS}. In contrast to \cite{PS}, we here do
not construct an approximate solution to the problem in $G^R$.
Instead, we use a simple estimate on its solutions which is uniform
with respect to $R \geq R_0$ and $\mu \in [\mu_1, \mu_2]$ and is
indifferent to the presence in $[\mu_1, \mu_2]$ of eigenvalues of
the operator $\{\mathcal{L} -\mu, \mathcal{B}\}$ in $G$.

The statement of basic boundary value problem, the precise
definition of scattering matrix, the detailed formulation of
computational method, and the statement of principal Theorem \ref
{PrinTh} are given in $\S 2$. The solvability of an auxiliary
problem in the domain $G^R$ is discussed in $\S 3$.  In $\S 4$, we
complete the justification of method  by proving Theorem
\ref{PrinTh}.
\section{Statement of the problem. Formulation of the
results}\label{s2}
\subsection{Boundary value problem}\label{ss2.1}
Let $G$ be a domain in $\mathbb{R}^{n+1}$ coinciding, outside a
large ball, with the union $\Pi_+^1\cup\ldots\cup\Pi_+^P$ of
non-overlapping semicylinders $\Pi_+^p=\{(y^p,t^p): y^p\in\Omega^p,
t^p>0\}$, where $(y^p,t^p)$ are local coordinates in  $\Pi_+^p$, the
cross-section $\Omega^p$ of $\Pi_+^p$ being a bounded domain in
$\mathbb{R}^{n}$. The boundary $\partial G$ of $G$ is supposed to be
smooth. We introduce a formally self-adjoint  $(k\times k)$-matrix
$\|\mathcal{L}_{ij}(x,D_x)\|$ of differential operators in $G$,
where $\mbox{ord}\,\mathcal{L}_{ij}=\tau_i+\tau_j$ with nonnegative
integers $\tau_j$ and $\tau_1+\ldots+\tau_k=m$. We also assume that
the Green formula
\begin{equation}\label{Green formula}
    (\mathcal{L}u,v)_G+(\mathcal{B}u,\mathcal{Q}v)_{\partial G}=(u,\mathcal{L}v)_G+(\mathcal{Q}u,\mathcal{B}v)_{\partial G}
\end{equation}
holds for all $u, v \in C_c^{\infty}(\overline{G})$, while the $(m
\times k)$-matrices $\mathcal{B}=\|\mathcal{B}_{qj}\|$ and
$\mathcal{Q}=\|\mathcal{Q}_{qj}\|$ consist of differential operators
such that $\mbox{ord}\,\mathcal{B}_{qj}=\sigma_q+\tau_j$, $\sigma_q$
being an integer,  and
$\mbox{ord}\,\mathcal{B}_{qj}+\mbox{ord}\,\mathcal{Q}_{qi}\leqslant
\tau_j+\tau_i-1$; here $(\cdot,\cdot)_G$ and
$(\cdot,\cdot)_{\partial G}$ stand for the inner products on
$L_2(G)$ and $L_2(\partial G)$. The coefficients of  $\mathcal{L}$,
$\mathcal{B}$, and $\mathcal{Q}$ are smooth in $\overline{G}$. We
suppose that the operator $\{\mathcal{L},\mathcal{B}\}$ of boundary
value problem in $G$ is elliptic.
\begin{rem}\label{R1}
{\rm We have used the (widest) ellipticity definition in the sense
of Agmon-Douglis-Nirenberg in order to include in consideration,
among others, some hydrodynamics problems. The other ellipticity
definitions can be obtained from that as special cases. For
instance, the scalar case corresponds to $k=1$, $\tau_1=m$, and
${\rm ord}\mathcal L =2m$; the boundary operator $\mathcal B$ is a
column $\{\mathcal B_1, \dots, \mathcal B_m \}$ with ${\rm ord}
\mathcal B_h =m_h$. A detailed description of various ellipticity
definitions as well as examples can be found, e.g.,  in \cite{A}}.
\end{rem}
Let us describe the coefficients of  $\mathcal{L}$ and
$\mathcal{B}$ in a neighborhood of infinity. Denote by
$\{L^p,B^p\}$ an elliptic boundary value problem operator in the
cylinder $\Pi^p=\Omega^p\times\mathbb{R}$ with coefficients
independent of $t^p\in\mathbb{R}$ and smooth in
$\overline{\Omega^p}$, while
$\mbox{ord}\,L^p_{ij}=\mbox{ord}\,\mathcal{L}_{ij}$ and
$\mbox{ord}\,B^p_{ij}=\mbox{ord}\,\mathcal{B}_{ij}$. We assume
that any coefficient $a$ of $\{\mathcal{L},\mathcal{B}\}$
satisfies
\begin{equation}\label{stabilization}
    D^{\alpha}\left(a(y^p,t^p)-a^p(y^p)\right)=O(\exp(-\delta t^p))
\end{equation}
in $\Pi^p_+$ for $t^p\rightarrow+\infty$, where $a^p$ is the
corresponding coefficient of  $\{L^p,B^p\}$ and $\delta$ is a
positive number. From ($\ref{Green formula}$) it follows that
there holds  the Green formula in every cylinder $\Pi^p$ obtained
from ($\ref{Green formula}$) by changing $G$ for $\Pi^p$ and
$\mathcal{L}$, $\mathcal{B}$, and $\mathcal{Q}$ for $L^p$, $B^p$,
and $Q^p$.

We consider the boundary value problem
\begin{equation}\label{Problem in G}
    \left\{%
\begin{array}{ll}
    \mathcal{L}(x,D_x)u(x)-\mu u(x)=0, & \hbox{$x\in G$,} \\
    \mathcal{B}(x,D_x)u(x)=0, & \hbox{$x\in\partial G$,} \\
\end{array}%
\right.
\end{equation}
with spectral parameter $\mu\in\mathbb{C}$. A number  $\mu$ is
called an eigenvalue of the operator $\{\mathcal{L},\mathcal{B}\}$,
if there exists a nonzero function $u\in L_2(G)$ smooth in
$\overline{G}$ and satisfying ($\ref{Problem in G}$). Such a
function is called an eigenfunction corresponding to the eigenvalue
$\mu$. From the Green formula it follows that any eigenvalue is
real. Every eigenfunction admits the estimate
$u(y^p,t^p)=O(\exp(-\varkappa t^p))$ for $t^p\rightarrow+\infty$ in
each $\Pi^p_+$ with a certain $\varkappa>0$. For any eigenvalue
there exist at most finitely many linearly independent
eigenfunctions.

\subsection{Space of waves. Scattering matrix}\label{ss2.2}

To simplify notation, we sometimes drop the superscript $p$, if the
context excludes misunderstanding. In every domain $\Omega=\Omega^p$
we introduce the operator pencil
$$\mathbb C \ni
\lambda\mapsto\mathfrak{A}(\lambda,\mu)=\{L(\lambda)-\mu
I,B(\lambda)\},
$$
where $L(\lambda)=L(y,D_y,\lambda)$, $B(\lambda)=B(y,D_y,\lambda)$,
so we have the pencils $\mathfrak{A}^1(\lambda,\mu)$, $\ldots$,
$\mathfrak{A}^P(\lambda,\mu)$. Let us fix, for the time being, the
parameter $\mu$. Considering $\lambda$ as spectral parameter of the
pencils, we shall use the same terminology as in $\cite{GGK}$.
 The spectrum of $\mathfrak{A}(\cdot,\mu)$ is symmetric about the real axis and
consists of normal eigenvalues, that is, isolated eigenvalues of
finite algebraic multiplicity. Every strip
$\{\lambda\in\mathbb{C}:|\mbox{Im}\,\lambda|<h<\infty\}$ contains at
most finitely many eigenvalues. It is known (see $\cite{NP}$) that
the total algebraic multiplicity of  the eigenvalues of
$\mathfrak{A}(\cdot, \mu)$ in the strip
$\{\lambda\in\mathbb{C}:|\mbox{Im}\,\lambda|<h\}$ is even for any
$h>0$; in particular, the total algebraic multiplicity of the real
eigenvalues is even.

If $\mu$ is not a threshold, then for any real eigenvalue of
$\mathfrak{A}^p(\cdot, \mu)$ there exist only eigenvectors and no
generalized eigenvectors.  We number all real eigenvalues of the
pencil $\mathfrak{A}^p(\cdot, \mu)$  counted according to their
(geometric) multiplicity. Let $\lambda_1^p, \ldots,
\lambda_{2M^p}^p$ be all such eigenvalues and let
$\varphi_1^p,\ldots,\varphi_{2M^p}^p$ be the corresponding
eigenvectors. The functions
$$
\Pi^p \ni (y, t) \mapsto
u_k^p(y,t)=\exp(i\lambda_k^pt)\varphi_k^p(y)
$$
satisfy the homogeneous problem
\begin{equation}\label{HomProblInPi}
\left\{%
\begin{array}{ll}
    (L^p(y,D_y,D_t)-\mu)v(y,t)=0, & \hbox{$(y,t)\in \Pi^p$;} \\
    B^p(y,D_y,D_t)v(y,t)=0, & \hbox{$(y,t)\in\partial \Pi^p$.} \\
\end{array}%
\right.
\end{equation}
Introduce
$$
q^p(u,v):=(L^pu,v)_{\Pi^p}+(B^pu,Q^pv)_{\partial\Pi^p}-(u,L^pv)_{\Pi^p}-(Q^pu,B^pv)_{\partial\Pi^p}.
$$
If $u,v\in C_c^{\infty}(\overline{\Pi^p})$, then $q^p(u,v)=0$.
Assume that $\chi\in C^{\infty}(\mathbb{R})$, $\chi(t)=1$ for
$t\geqslant 2$ and $\chi(t)=0$ for $t\leqslant 1$. The form $q^p$
extends to the pairs $\{\chi u_j^p,\chi u_k^p\}$, and there is valid
the following assertion (see \cite{NP}).
\begin{prop}
The eigenvectors $\{\varphi_k^p\}$ can be chosen to satisfy
 $q^p(\chi u_j^p,\chi u_k^p)=\pm i\delta_{jk}$,
while the sign is determined and cannot be taken arbitrarily.
\end{prop}
We extend the functions  $\chi u_k$ by zero from  $\Pi^p_+$ to the
domain $G$ keeping the same notation for the extended functions.
Introduce a space $E$ of $w\in C^{\infty}(\overline{G})$
satisfying $D^{\alpha}w(x)=O(\exp(-\beta|x|))$ as
$|x|\rightarrow+\infty$ for any multiindex $\alpha$ and some
positive $\beta<\delta$, where $\delta$ is the number in
(\ref{stabilization}). Let $\mathcal{W}$ stand for the linear span
of functions $\chi u_k+w$ with $w\in E$. The form
\begin{equation}\label{q}
q(u,v)=(\mathcal{L}u,v)_G+(\mathcal{B}u,\mathcal{Q}v)_{\partial
G}-(u,\mathcal{L}v)_G-(\mathcal{Q}u,\mathcal{B}v)_{\partial G}
\end{equation}
takes finite values for any $u,v\in\mathcal{W}$. It is evident
that $q(u,v)=-\overline{q(v,u)}$ and $q(u,u)\in i\mathbb{R}$. By
definition, $u\in\mathcal{W}$ is an incoming (outgoing) wave, if
$iq(u,u)>0$ ($iq(u,u)<0$). In the space $\mathcal{W}$ one can
choose the basis
\begin{equation}\label{basis of waves}
    u_1^+,\ldots,u_M^+,u_1^-,\ldots,u_M^-,\qquad M=\sum_{p=1}^PM^p,
\end{equation}
subject to the conditions
\begin{equation}\label{orthogonality of waves}
    q(u_j^{\pm},u_k^{\pm})=\mp i\delta_{jk},\quad
    q(u_j^{\pm},u_k^{\mp})=0,\quad j,k=1,\ldots,M;
\end{equation}
here $u_1^+,\ldots,u_M^+$ are incoming waves and
$u_1^-,\ldots,u_M^-$ are outgoing ones.

Let $\gamma$ be a positive number such that $\gamma <\delta$ and the
strip $\{\lambda\in\mathbb{C}:|\mbox{Im}\,\lambda|\leq \gamma) \}$
contains only real eigenvalues of the pencils
$\mathfrak{A}^p(\cdot,\mu)$, $p=1, \dots, P$,  for all $\mu \in
[\mu_1, \mu_2]$. It is known \cite{NP} that in the space of bounded
solutions to the problem (\ref{Problem in G}) there exist elements
$Y_1 (\cdot, \mu) \ldots,Y_M (\cdot, \mu)$ such that
\begin{equation}\label{Y}
    Y_j(x, \mu)=u_j^+(x, \mu)+\sum_{k=1}^MS_{jk}(\mu)u_k^-(x, \mu)+O(e^{-\gamma |x|})
\end{equation}
as $|x|\rightarrow\infty$. If the number $\mu$ is not an eigenvalue
of the problem (\ref{Problem in G}), then $Y_j (\cdot, \mu)$ are
uniquely determined and form a basis in the space of bounded
solutions to the homogeneous problem (\ref{Problem in G}), that is,
in the space of continuous spectrum eigenfunctions. Otherwise, any
$Y_j (\cdot, \mu)$ is determined up to an eigenfunction of
(\ref{Problem in G}) belonging to $L_2(G)$ . Then any bounded
solution of (\ref{Problem in G}) can be represented by a linear
combination of $Y_j (\cdot, \mu)$ up to an eigenfunction in $L_2
(G)$.

The matrix $S(\mu)=\|S_{jk}(\mu)\|_{j,k=1}^M$  in (\ref{Y}) has
uniquely been determined for all $\mu \in [\mu_1, \mu_2]$; it is
independent of the arbitrariness in the definition of $Y_j (\cdot,
\mu)$  when $\mu$ is an eigenvalue. The matrix  $S(\mu)$ is called
the scattering matrix. It is unitary for all  $\mu$.

In what follows we do not usually indicate the dependence of
$u_j^{\pm}$,  $Y_j$, etc., on the spectral parameter $\mu$. The
context excludes misunderstanding.

\subsection{Method for computing the scattering matrix. Formulation of the principal theorem}\label{ss2.3}

For large $R$ we introduce
$$
\Pi_+^{p,R}=\{(y^p,t^p)\in\Pi^p:t^p>R\},\quad
G^R=G\setminus\cup_{p=1}^N\Pi_+^{p,R}
$$
and set $$\partial G^R\setminus\partial
G=\Gamma^R=\cup_p\Gamma^{p,R}, \quad
\Gamma^{p,R}=\{(y^p,t^p)\in\Pi^p:t^p=R\}.$$ Then
\begin{equation}\label{Green formula in GR}
    (\mathcal{L}u,v)_{G^R}+(\mathcal{B}u,\mathcal{Q}v)_{\partial
    G^R\setminus\Gamma^R}+(\mathcal{N}u,\mathcal{D}v)_{\Gamma^R}
    =(u,\mathcal{L}v)_{G^R}+(\mathcal{Q}u,\mathcal{B}v)_{\partial
    G^R\setminus\Gamma^R}+(\mathcal{D}u,\mathcal{N}v)_{\Gamma^R},
\end{equation}
where $\mathcal{D}$ and $\mathcal{N}$ are $(m\times k)$--matrices
of differential operators, $\mathcal{D}$ being a Dirichlet system
(see \cite{LM}). An example of Dirichlet system is presented by
the matrix consisting of $m$ rows of the form
$e^{(j)}\partial_{\nu}^h$, where $j=1,\ldots,k$,
$h=1,\ldots,\tau_j-1$,
$e^{(j)}=(\delta_{1,j},\ldots,\delta_{k,j})$, and  $\nu$ is the
outward normal to $\Gamma^R$.

We look for the row  $(S_{l1},\ldots,S_{lM})$ of the scattering
matrix $S=S(\mu)$. As approximation to the row, we take the
minimizer of a quadratic functional. To construct such a functional
we consider the problem
\begin{eqnarray}\label{problem in GR for X}\nonumber
  (\mathcal{L}(x,D_x)-\mu)\mathcal{X}_l^R &=& 0,\quad x\in G^R, \\
  \mathcal{B}(x,D_x)\mathcal{X}_l^R &=& 0,\quad x\in \partial
    G^R\setminus\Gamma^R,\\\nonumber
  (\mathcal{N}+i\zeta\mathcal{D})\mathcal{X}_l^R &=&
  (\mathcal{N}+i\zeta\mathcal{D})(u_l^++\sum\nolimits_{j=1}^Ma_ju_j^-),\,x\in\Gamma^R,
\end{eqnarray}
where $\zeta$ is a fixed number in $\mathbb{R}\setminus\{0\}$ and
$a_1,\ldots,a_M$ are complex numbers.

Let us explain the origin of the problem. A solution  $Y_l$ of the
homogeneous problem (\ref{Problem in G}) satisfies the two first
equations (\ref{problem in GR for X}). The asymptotics (\ref{Y}) can
be differentiated so
$$
(\mathcal{N}+i\zeta\mathcal{D})Y_l=(\mathcal{N}+i\zeta\mathcal{D})(u_l^++\sum_{j=1}^Ma_ju_j^-)
+O(e^{-\gamma R})
$$
for $a_j=S_{lj}$. Thus $Y_l$ gives an exponentially small
discrepancy to the last equation (\ref{problem in GR for X}). As
approximation to the row $(S_{l1},\ldots,S_{lM})$, we take the
minimizer $a^0(R)=(a^0_1(R),\ldots,a^0_M(R))$ of the functional
\begin{equation}\label{functional}
(a_1,\ldots,a_M)\mapsto
J_l^R(a_1,\ldots,a_M)=\|\mathcal{D}(\mathcal{X}_l^R-u_l^+-\sum_{j=1}^Ma_ju_j^-);L_2(\Gamma^R)\|^2,
\end{equation}
where $\mathcal{X}_l^R$ is a solution to the problem (\ref{problem
in GR for X}). One can expect that  $a_j^0(R)\rightarrow S_{lj}$
with exponential rate as $R\rightarrow\infty$ and $j=1,\ldots,M$.

To clarify the dependence of $\mathcal{X}_l^R$ on the parameters
$a_1,\ldots,a_M$, we consider the problems
\begin{eqnarray}\label{problem in GR for v+-}\nonumber
  (\mathcal{L}(x,D_x)-\mu)v_j^{\pm} &=& 0,\quad x\in G^R; \\
  \mathcal{B}(x,D_x)v_j^{\pm} &=& 0,\quad x\in \partial
    G^R\setminus\Gamma^R;\\\nonumber
  (\mathcal{N}+i\zeta\mathcal{D})v_j^{\pm} &=&
  (\mathcal{N}+i\zeta\mathcal{D})u_j^{\pm},\,x\in\Gamma^R;\quad
  j=1,\ldots,M.
\end{eqnarray}
It is evident that $\mathcal{X}_l^R=v_{l,R}^++\sum_ja_jv_{j,R}^-$,
where $v_j^{\pm}=v_{j,R}^{\pm}$ are solutions to (\ref{problem in GR
for v+-}). We introduce the $(M\times M)$--matrices with entries
\begin{equation}\label{EF}
    \begin{split}
    \mathcal{E}_{ij}^R=\left(\mathcal{D}(v_i^--u_i^-),\mathcal{D}(v_j^--u_j^-)\right)_{\Gamma^R},\\
    \mathcal{F}_{ij}^R=\left(\mathcal{D}(v_i^+-u_i^+),\mathcal{D}(v_j^--u_j^-)\right)_{\Gamma^R}
\end{split}
\end{equation}
and set
$$
\mathcal{G}_i^R=\left(\mathcal{D}(v_i^+-u_i^+),\mathcal{D}(v_i^+-u_i^+)\right)_{\Gamma^R}.
$$
Now the functional (\ref{functional}) can be written in the form
$$
J_l^R(a)=\langle
a\mathcal{E}^R,a\rangle+2\mbox{Re}\,\langle\mathcal{F}_l^R,a\rangle+\mathcal{G}_l^R,
$$
where $\mathcal{F}_l^R$ is the $l$-th row of the matrix
$\mathcal{F}^R$ and  $\langle\cdot,\cdot\rangle$ is the inner
product on $\mathbb{C}^M$. The minimizer $a^0$ (a row) satisfies
$a^0(R)\mathcal{E}^R+\mathcal{F}^R_l=0$. Therefore, as approximation
$S^R(\mu)$ to the scattering matrix $S(\mu)$ we have a solution of
the equation $S^R\mathcal{E}^R+\mathcal{F}^R=0$.

\begin{thm}\label{PrinTh}
Let $\zeta $ be any fixed number in $\mathbb R\setminus \{0\}$ and
let the interval $[\mu_1, \mu_2]$ of the continuous spectrum of
problem (\ref {Problem in G}) be free from the threshold values of
the spectral parameter $\mu$. Then for all $\mu \in [\mu_1, \mu_2]$
and $R>R_0$, where $R_0$ is a sufficiently large number, there
exists a unique minimizer $a^0(R, \mu)=(a_1^0(R, \mu),\ldots,$
$a_M^0(R, \mu))$ of the functional $J_l^R (\cdot; \mu)$ in
(\ref{functional}). The estimates
$$
|a_j^0(R, \mu)-S_{lj}(\mu)|\leqslant Ce^{-\Lambda R},\quad
j=1,\ldots,M,
$$
hold with constant $C$ independent of $\mu$ and $R$, while
$\Lambda=\min\{\beta,\gamma\}$, $\beta$ is a positive number less
than $\delta$ in (\ref{stabilization}), $\gamma$ is the same as in
(\ref{Y}).
\end{thm}

\section{Problem in the domain $G^R$}\label{s3.1}

Introduce the boundary value problem
\begin{eqnarray}\label{3.1}
\mathcal{L} (x, D_x)u(x)-\mu u(x)&=f(x),&  \,\, x\in G^R,  \nonumber \\
\mathcal{B} (x, D_x) u(x)&=g(x),& \, \, x\in \partial G^R\setminus
\Gamma^R,
\nonumber \\
(\mathcal{N} (x, D_x) + i\zeta \mathcal{D} (x, D_x))u(x)&=h(x),&
\,\, x\in \Gamma^R,
\end{eqnarray}
where $\zeta \in \mathbb{R} \setminus \{0\}$ and $\mu \in
\mathbb{R}$. In this section we discuss the unique solvability of
the problem.

The boundary $\partial G^R$ contains  $\partial \Gamma^R$, which is
an edge for $\dim G>2$ or the union of corner points for $\dim G=2$.
The boundary conditions have discontinuities along $\partial
\Gamma^R$. When studying problem (\ref{3.1}), one can use the
traditional scheme of the theory of elliptic boundary value problems
in domains with piecewise smooth boundary (see, e.g., \cite{NP},
\cite{Gris}, \cite{KMR}). In contrast to the smooth situation, the
choice of function spaces for a boundary value problem has not been
universal and requires taking into account the specific properties
of solutions near the edges; sometimes weighted function spaces turn
out to be suitable for the purpose. For the classical problems in
mathematical physics the needed spaces are known. For those reasons
we restrict ourselves to considering some specific features of
problem (\ref{3.1}), examples, and postulating the needed properties
of function spaces to be used.

Assume that the function spaces for problem (\ref{3.1}) have been
chosen so that the operator of the problem  is Fredholm (having
closed range and finite-dimensional kernel and cokernel). Then the
triviality of kernel and cokernel is necessary and sufficient for
the existence of a unique solution to the problem with any
right-hand side. To prove the triviality, note that  (\ref{Green
formula in GR}) leads to
\begin{eqnarray}\label{3.1a}
(\mathcal{L} -\mu)u, v)_{G^R}+(\mathcal{B} u, \mathcal{Q}
v)_{\partial G^R \setminus \Gamma ^R} + ((\mathcal{N} + i\zeta
\mathcal{D})u, \mathcal{D} v)_{\Gamma ^R}\nonumber
\\=(u, (\mathcal{L} -\mu) v)_{G^R} + (\mathcal{Q} u, \mathcal{B} v)_{\partial G^R \setminus
\Gamma ^R} + (\mathcal{D} u, (\mathcal{N} -i\zeta
\mathcal{D})v)_{\Gamma^R}.
\end{eqnarray}
Let $w$ be a solution to the homogeneous problem (\ref{3.1}).
Setting $u=v=w$ in (\ref{3.1a}), we obtain
\begin{equation}\label{3.1b}
\|\mathcal{D} w; L_2 (\Gamma^R)\|=0.
\end{equation}
This and the homogeneous boundary condition  $(\mathcal{N} (x, D_x)
+ i\zeta \mathcal{D} (x, D_x))w(x)=0$ for  $x \in \Gamma^R$ imply
that $w$ has zero Cauchy data at $\Gamma^R$. If the coefficients of
the operator are sufficiently smooth for the applicability of the
unique continuation theorem (see \cite{BJS}, part II, \S 1.4), then
we obtain the triviality of kernel. Similar considerations for the
adjoint problem provide the triviality of cokernel. It is supposed
that the function spaces where the problem (\ref{3.1}) is considered
admit such a reasoning. We illustrate this scheme of the analysis of
problem (\ref{3.1}) by the two following examples.

\noindent {\bf Example 1.} {\rm We assume that $\dim G=2$ and
consider the problem}
\begin{eqnarray}\label{3.1c}%{3.10}
(\Delta -\mu)u(x)&=& f(x), \,\,\, x\in G^R, \nonumber \\
u(x)&=& g(x), \,\,\, x\in \partial G^R \setminus \Gamma^R, \nonumber \\
\partial_\nu u(x)+i\zeta u(x)&=&h(x), \,\,\, x\in \Gamma^R.
\end{eqnarray}
{\rm With every corner point of the boundary $\partial G^R$ we
associate the problem with complex parameter (operator pencil)}
\begin{eqnarray}\label{3.1d} %{3.11}
(\partial_\omega^2 -\lambda^2)v (\omega)&=&0, \,\,\, \omega
\in (0, \pi/2),  \nonumber \\
v(0)= v'(\pi/2)&=&0.
\end{eqnarray}
The spectrum of problem (\ref{3.1d}) consists of simple eigenvalues
$\lambda_q=(2q+1)i$, where $q=0, \pm 1, \dots$, while $\omega
\mapsto \sin (2q+1)\omega$ is an eigenfunction corresponding to
$\lambda_q$.

Introduce the space $V^l_\beta (G^R)$ with norm
$$
\|u; V^l_\beta (G^R)\|=\left( \sum_{|\alpha|\leq l} \int_{G^R}
r^{2(\beta -l+|\alpha|)}|D_x^\alpha u(x)|^2\,dx \right)^{1/2},
$$
where $\beta \in \mathbb{R}$, $l=0, 1, \dots$, and $r$ denotes a
function that coincides near a corner point with the distance to the
point, equals 1 outside a neighborhood of the corner points being
smooth and positive  on $\overline {G^R}$ (except at the corner
points). Let also  $V^{l-1/2}_\beta (\partial G^R\setminus
\Gamma^R)$ and $V^{l-1/2}_\beta (\Gamma^R)$ with  $l=1, 2, \dots$
stand for the space of traces of the functions in $V^{l}_\beta
(G^R)$ on $\partial G^R\setminus \Gamma^R$ and $\Gamma^R$
respectively.

The operator $\mathcal{A}^R(\mu)$ of problem (\ref{3.1}) implements
a continuous mapping
\begin{eqnarray}\label{3.1e} %{3.12}
V^2_\beta (G^R) \ni u \mapsto \mathcal{A}^R(\mu)u=\{f, g, h\}\in
V^0_\beta (G^R)\times V^{3/2}_\beta (\partial G^R \setminus
\Gamma^R)\times V^{1/2}_\beta (\Gamma^R).
\end{eqnarray}
It is known that the operator (\ref{3.1e}) is Fredholm if and only
if $\beta -1$ coincides with none of the numbers ${\rm Im}\,
\lambda_q$, that is, $\mathcal{A}^R(\mu)$ is Fredholm when $\beta$
is not even. If $w$ satisfies the homogeneous problem (\ref{3.1})
and  $w\in V^2_\beta (G^R)$ with a certain $\beta \in (2q, 2q+2)$
for an integer $q$, then near a corner point
$$w (x) =Cr^{2q+1}\sin(2q+1)\omega + O(r^{2q+2-\varepsilon}),$$
where $r, \omega$ are polar coordinates centered at the corner
point, $C$ is a constant, and $\varepsilon$ is any positive number
subject to $\varepsilon <1$. Hence for every element in the kernel
of the operator  (\ref{3.1e}) with $\beta \in (0, 2)$, the formula
(\ref{3.1b}) holds, which now takes the form $\|w; L_2
(\Gamma^R)\|=0$. Therefore the kernel is trivial for $\beta \in (0,
2)$ and consequently for all $\beta <2$.

We turn to the cokernel. Denote by  $V^{-l}_{-\beta} (G^R)$ the
space adjoint to $V^l_\beta (G^R)$ with respect to the inner
product on $L_2(G^R)$. Let  $V^{-l-1/2}_{-\beta}(\Gamma^R)$ stand
for the space adjoint to  $V^{l+1/2}_\beta (\Gamma^R)$ with
respect to the inner product on $L_2 (\Gamma^R)$, $l=0, 1$.
Finally, denote by
 $\mathcal{A}^R(\mu)^*$ the operator adjoint to the operator
(\ref{3.1e}),
\begin{equation}\label{3.1f}%{3.13}
\mathcal{A}^R(\mu)^*: V^0_{-\beta} (G^R)\times
V^{-3/2}_{-\beta}(\partial G^R \setminus \Gamma^R) \times
V^{-1/2}_{-\beta}(\Gamma^R) \to V^{-2}_{-\beta}(G^R).
\end{equation}
The cokernel of operator (\ref{3.1e}) coincides with the kernel of
operator (\ref{3.1f}). According to the known results on the
regularity of elliptic problem solutions, for any element   $\{ u,
v, w\}$ in the kernel of (\ref{3.1f}) there is the inclusion
\begin{equation}\label{3.1g}%{3.14}
\{u, v, w\} \in V^2_{2-\beta} (G^R)\times V^{1/2}_{2-\beta}(\partial
G^R \setminus \Gamma^R) \times V^{3/2}_{2-\beta}(\Gamma^R),
\end{equation}
the function $u$ satisfies the homogeneous problem (\ref{3.1c}) with
$\partial_\nu +i\zeta$ replaced for $\partial_\nu -i\zeta$ in the
boundary condition at $\Gamma^R$, while $v$ and $w$ are determined
by
\begin{eqnarray}\label{3.1h}%{3.15}
v(x)=-\partial_\nu u(x), \,\, x\in \partial G^R\setminus \Gamma^R;
\,\,\, w(x)=u(x), \,\, x\in \Gamma^R.
\end{eqnarray}
The above discussion of the kernel triviality of operator
 (\ref{3.1e}) does not depend on a sign of
$\zeta$. Therefore taking into account (\ref{3.1g}), we obtain $u=0$
for all $\beta$ such that  $2-\beta <2$. By virtue of  (\ref{3.1h}),
for the same $\beta$ we have $v=0$ and $w=0$. Thus if $\beta \in (0,
2)$, then both kernel and cokernel of (\ref{3.1e}) are trivial. It
follows that for $\beta \in (0, 2)$ and for all $\mu$ and $\zeta
\neq 0$, the operator (\ref{3.1e}) is an isomorphism.
 Moreover, it can be shown that, for the even numbers $\beta$, the range of the operator is not closed,
 the   cokernel is nontrivial for $\beta <0$, and
 the  kernel is nontrivial for  $\beta >2$.  $\Box$

Slightly modifying the statement of problem (\ref{3.1}), it is
sometimes possible to do without edges at the boundary and
discontinuities in the boundary conditions. Then the analysis of the
problem becomes simpler, while all the rest in the sequel requires
no essential changes.

\noindent {\bf Example 2.} Let the role of initial problem
(\ref{Problem in G}) be played by the Neumann problem
\begin{eqnarray*}
(\mathcal{L} (x, D_x) - \mu) u(x)&=&f(x),\,\,\, x \in G, \\
\mathcal{N} (x, D_x)u(x)&=&g(x),\,\,\, x \in \partial G,
\end{eqnarray*}
where $\mathcal{L}$ is the same matrix differential operator as in
(\ref{Problem in G}). Denote by $\widetilde{G}^R$ the bounded domain
with smooth boundary obtained from $G$ by cutting off the
cylindrical end $\Pi^p_+$ by a surface $\widetilde {\Gamma}^{p,R}$
such that $\widetilde \Gamma^{p, R} \subset \{(y^p, t^p)\in \Pi^p_+:
R<t^p<R+1\}$, $p=1, \dots, P$. As $R$ varies,  the surface
$\widetilde {\Gamma}^{p,R}$ moves in a parallel way along the axis
of $\Pi^p_+$. We set $\widetilde{\Gamma}^R = \cup_p
\widetilde{\Gamma}^{p, R}$. Let us choose a smooth cut-off function
$\chi$ on  $\partial \widetilde{G}^R$ such that ${\rm supp}\chi
\subset \widetilde {\Gamma}^R$ and ${\rm mes}\{x \in
\widetilde{\Gamma}^R : \chi (x)=1\} >0$. Instead of (\ref{3.1}), we
introduce the problem
\begin{eqnarray}\label{3.1i}
\mathcal{L} (x, D_x)u(x)-\mu u(x)&=f(x),&  \,\, x\in \tilde{G}^R,  \nonumber \\
(\mathcal{N} (x, D_x) + i\zeta \chi \mathcal{D} (x,
D_x))u(x)&=g(x),& \,\, x\in
\partial \tilde{G}^R.
\end{eqnarray}
There is the Green formula
\begin{eqnarray}\label{3.1j}%{3.17}
((\mathcal{L} -\mu) u , v)_{\tilde{G}^R}+((\mathcal{N} +i\zeta \chi
\mathcal{D})u, \mathcal{D} v)_{\partial \tilde{G}^R}\nonumber \\ =
(u, (\mathcal{L} -\mu)v)_{\tilde{G}^R}+(\mathcal{D} u, (\mathcal{N}
-i\zeta \chi \mathcal{D})v)_{\partial \tilde{G}^R}.
\end{eqnarray}
The operator $\mathcal{A}^R(\mu)$ of problem (\ref{3.1i}) implements
a continuous mapping
\begin{eqnarray}\label{3.1k}%{3.18}
\mathcal{A}^R(\mu): \prod_{j=1}^k H^{l+\tau+\tau_j}(\widetilde{G}^R)
\to \prod_{j=1}^k H^{l+\tau -\tau_j}(\widetilde{G}^R)\times
\prod_{h=1}^m H^{l-\sigma_h-1/2}(\partial \widetilde{G}^R),
\end{eqnarray}
where $H^s(\widetilde{G}^R)$ and $H^s(\partial \widetilde{G}^R)$ are
the usual Sobolev spaces, $\tau =\max\{\tau_1, \dots, \tau_k\}$,
$l\geq \max\{1+\max\sigma_h, 0\}$, and the numbers $k$ and $m$ were
defined for the matrix $\mathcal L (x, D_x)$ in Section \ref{ss2.1};
in the scalar case, $k=1$, $\tau =\tau_1=m$, and $\sigma_h=m_h-m$
(see \ref{ss2.1}). The operator (\ref{3.1k}) is Fredholm. From
(\ref{3.1j}) it follows that a solution $u$ of the homogeneous
problem (\ref{3.1i}) satisfies $\chi (x)\mathcal{D} (x, D_x) u(x)=0$
on $\partial \tilde{G}^R$. This and the homogeneous boundary
condition in (\ref{3.1i}) imply that $u$ has zero Cauchy data on the
set $\{x\in
\partial \tilde{G}^R: \chi (x)=1\}$. Under the condition of unique continuation
theorem, it follows that the kernel of operator (\ref{3.1k}) is
trivial. The triviality of cokernel can be proved in a similar way.

As $\{\mathcal{L}, \mathcal{B}\}$, one can take, for instance, the
operator $\{ \mathcal{E}, \mathcal{N} \}$ of elasticity theory; here
$$
\mathcal{E} (x, D_x)u(x)=-\mu \nabla_x\cdot \nabla_x
-(\lambda+\mu)\nabla_x\cdot u(x)
$$
is the Lam\'{e} system with parameters $\lambda$ and $\mu$,
$\mathcal{N} (x, D_x)u(x)=\sigma (u; x)\nu (x)$, $\nu$ is the
outward normal, and $\sigma (u; x)$ is the stress tensor,
$$
\sigma_{jk}(u, x)=\mu \left( \frac{\partial u_j}{\partial x_k} +
\frac{\partial u_k}{\partial x_j} \right) +\lambda
\delta_{jk}\nabla_x\cdot u. \,\,\,\Box
$$

\section{Justification of the method for computing scattering matrices}\label{s4}

To justify the method, we have to verify that the matrix
$\mathcal{E}^R$ with entries (\ref{EF}) is nonsingular and the
minimizer $a^0(R)$ of (\ref{functional}) tends to the $l$-th row of
scattering matrix as $R\rightarrow\infty$.
\begin{prop}\label{p4.1}
Let $u_j^{\pm}$ be incoming and outgoing waves (\ref{basis of
waves}) satisfying (\ref{orthogonality of waves}). Then
\begin{eqnarray*}
  &(\mathcal{N}u_j^{\pm},
\mathcal{D}u_k^{\pm})_{\Gamma^{R}}- (\mathcal{D}u_j^{\pm},
\mathcal{N}u_k^{\pm})_{\Gamma^{R}}=\mp i\delta_{jk}+O(e^{-\beta R}), \\
  &(\mathcal{N}u_j^{\pm},
\mathcal{D}u_k^{\mp})_{\Gamma^{R}}- (\mathcal{D}u_j^{\pm},
\mathcal{N}u_k^{\mp})_{\Gamma^{R}}=O(e^{-\beta R})
\end{eqnarray*}
for $R\rightarrow+\infty$, where $\beta$ is the number in the
definition of $\mathcal{W}$ in 2.1, i. e., $0<\beta<\delta$, where
$\delta$ is in (\ref{stabilization}).
\end{prop}
\noindent {\bf Proof.} According to  (\ref{Green formula in GR}),
$$
(\mathcal{N}u,\mathcal{D}v)_{\Gamma^R}-(\mathcal{D}u,\mathcal{N}v)_{\Gamma^R}
    =(u,\mathcal{L}v)_{G^R}+(\mathcal{Q}u,\mathcal{B}v)_{\partial
    G^R\setminus\Gamma^R}-(\mathcal{L}u,v)_{G^R}-(\mathcal{B}u,\mathcal{Q}v)_{\partial
    G^R\setminus\Gamma^R}.
$$
If $u$ and $v$ are elements of the lineal $\mathcal{W}$ introduced
in \ref{ss2.2}, then the right-hand side of the last equality
differs from $q(u,v)$ by a term $O(e^{-\beta R})$ as
$R\rightarrow+\infty$. Hence
\begin{eqnarray*}
  &(\mathcal{N}u_j^{\pm},
\mathcal{D}u_k^{\pm})_{\Gamma^{R}}- (\mathcal{D}u_j^{\pm},
\mathcal{N}u_k^{\pm})_{\Gamma^{R}}=q(u_j^{\pm},u_k^{\pm})+O(e^{-\beta R}), \\
  &(\mathcal{N}u_j^{\pm},
\mathcal{D}u_k^{\mp})_{\Gamma^{R}}- (\mathcal{D}u_j^{\pm},
\mathcal{N}u_k^{\mp})_{\Gamma^{R}}=q(u_j^{\pm},u_k^{\mp})+O(e^{-\beta
R}).
\end{eqnarray*}
It remains to take into account (\ref{orthogonality of waves}).
$\Box$

Let $u$ be a solution to a problem of the form $\mathcal
{A}^R(\mu)w=(0, 0, h)$ with $(0, 0, h)\in \mathcal {R}(\mathcal{A}^R
(\mu) \cap L_2 (\Gamma^R)$, where $\mathcal {R}(\mathcal{A}^R
(\mu))$ stands for the range of the operator of problem (\ref{3.1}).
In what follows, we assume that $\mathcal D u\in L_2 (\Gamma^R)$.
Moreover, we also assume that if $v$ possesses similar properties,
then for $u$ and $v$ there holds the Green formula (\ref{Green
formula in GR}). (Note that such assumptions have been fulfilled for
the problems in Examples 1 and 2, Section \ref{s3.1}.)
\begin{prop}\label{prop4.2}
The matrix $\mathcal{E}^R$ with entries (\ref{EF})  is nonsingular
for all  $R\geqslant R_0$, where  $R_0$ is a sufficiently large
number.
\end{prop}
\noindent {\bf Proof.} Suppose the proposition is false. Then for
any $R^0$ there exists a number $R>R^0$ such that the matrix
$\mathcal{E}^R$ is singular, while $\mathcal{U}=\sum_jc_ju_j^-$ and
$\mathcal{V}=\sum_jc_jv_j^-$ satisfy
\begin{equation}\label{5.1}
    \mathcal{D}\mathcal{U}=\mathcal{D}\mathcal{V}\quad\mbox{on}\;\Gamma^R,
\end{equation}
where $v_j^-$ is a solution to problem (\ref{problem in GR for v+-})
and $\overrightarrow{c}=(c_1,\ldots,c_M)$ with
$|\overrightarrow{c}|=1$. According to the equation on  $\Gamma^R$
in (\ref{problem in GR for v+-}), we have
\begin{equation}\label{5.2}
    \mathcal{N}\mathcal{U}=\mathcal{N}\mathcal{V}\quad\mbox{on}\;\Gamma^R.
\end{equation}
We set $u=v=\mathcal{V}$ in (\ref{Green formula in GR}), take
account of (\ref{5.1}), (\ref{5.2}) and of the two first equations
(\ref{problem in GR for v+-}) and obtain
\begin{equation}\label{5.3}
    (\mathcal{N}\mathcal{U},\mathcal{D}\mathcal{U})_{\Gamma^R}-
    (\mathcal{D}\mathcal{U},\mathcal{N}\mathcal{U})_{\Gamma^R}=0.
\end{equation}
This and Proposition \ref{p4.1} imply that
$$
0=i\sum_j|c_j|^2+o(1)=i+o(1),
$$
a contradiction. $\Box$
\begin{prop}\label{p4.3}
Let $u$ be a solution to the problem (\ref{3.1}) with right-hand
side $(0, 0, h)$,  while $h \in L_2(\Gamma^R)$. Then
\begin{equation}\label{est}
\|\mathcal{D}u;L_2(\Gamma^R)\|\leqslant\frac{1}{|\zeta|}\|h;L_2(\Gamma^R)\|.
\end{equation}
\end{prop}
\noindent {\bf Proof.} From (\ref{Green formula in GR}) it follows
that
\begin{eqnarray}\nonumber
    (\mathcal{L}u,v)_{G^R}+(\mathcal{B}u,\mathcal{Q}v)_{\partial
    G^R\setminus\Gamma^R}-(u,\mathcal{L}v)_{G^R}-(\mathcal{Q}u,\mathcal{B}v)_{\partial
    G^R\setminus\Gamma^R}=\\\nonumber
    =((\mathcal{N}+i\zeta\mathcal{D})u,\mathcal{D}v)_{\Gamma^R}-(\mathcal{D}u,(\mathcal{N}+
    i\zeta\mathcal{D})v)_{\Gamma^R}-2i\zeta(\mathcal{D}u,\mathcal{D}v)_{\Gamma^R}.
\end{eqnarray}
We set $v=u$ and obtain
$$
0=(h,\mathcal{D}u)_{\Gamma^R}-(\mathcal{D}u,h)_{\Gamma^R}-2i\zeta\|\mathcal{D}u;L_2(\Gamma^R)\|^2.
$$
Then
$$
2|\zeta|\|\mathcal{D}u;L_2(\Gamma^R)\|^2=|(h,\mathcal{D}u)_{\Gamma^R}-(\mathcal{D}u,h)_{\Gamma^R}|\leqslant
2\|\mathcal{D}u;L_2(\Gamma^R)\|\|h;L_2(\Gamma^R)\|
$$
and we arrive at (\ref{est}). $\Box$
\begin{prop}\label{p4.4}
Let  $a(R)=(a_1(R),\ldots,a_M(R))$ be a minimizer of  $J_l^R$ in
(\ref{functional}). Then
\begin{equation}\label{5.5}
    J_l^R\left(a(R)\right)=O(e^{-2\gamma
    R})\quad\mbox{as}\;R\rightarrow\infty,
\end{equation}
where $\gamma$ is the same as in (\ref{Y}). For all  $R\geqslant
R_0$,
$$
|a_j(R)|\leqslant \mbox{const}<\infty,\quad j=1,\ldots,M.
$$
\end{prop}
\noindent {\bf Proof.} Denote by  $Y_l^R$ a solution to the problem
(\ref{problem in GR for X}) with $a_j$, $j=1,\ldots,M$, equal to the
entries $S_{lj}$ of the scattering matrix  $S$ of problem
(\ref{Problem in G}). Since the asymptotics (\ref{Y}) can be
differentiated, we obtain
$$
(\partial_{\nu}+i\zeta)(Y_l^R-Y_l)|_{\Gamma}=O(e^{-\gamma R}).
$$
The difference $Y_l^R-Y_l$ satisfies the two first equations of the
problem (\ref{3.1}) with  $f=0$ and $g=0$, therefore  (\ref{est})
holds for $u=Y_l^R-Y_l$:
\begin{eqnarray*}
  && \|\mathcal{D}(Y_l^R-Y_l);L_2(\Gamma^R)\| \leqslant |\zeta|^{-1}\|(\mathcal{N}+i\zeta\mathcal{D})(Y_l^R-Y_l);L_2(\Gamma^R)\|
  \leqslant ce^{-\gamma R}.
\end{eqnarray*}
This and  (\ref{Y}) lead to the estimate
$$
J_l^R(S_l)=\|\mathcal{D}(Y_l^R-(u^+_l+\sum_{j=1}^MS_{lj}u_j^-));L_2(\Gamma^R)\|^2\leqslant
c e^{-2\gamma R}
$$
with constant  $c$ independent of  $R$. Owing to
$J_l^R(a(R))\leqslant J_l^R(S_l)$, we have (\ref{5.5}).

Let us estimate the minimizer $a(R)$. Denote by $Z_l^R$ the solution
of problem (\ref{problem in GR for X}) corresponding to
$a(R)=(a_1(R),\ldots,a_M(R))$. We set $u=v=Z_l^R$ in (\ref{Green
formula in GR}) and obtain
\begin{equation}\label{5.6}
    (\mathcal{N}Z_l^R,\mathcal{D}Z_l^R)_{\Gamma^R}-
    (\mathcal{D}Z_l^R,\mathcal{N}Z_l^R)_{\Gamma^R}=0.
\end{equation}
By virtue of  (\ref{5.5})
\begin{equation}\label{5.7}
    \|\mathcal{D}(Z_l^R-(u^+_l+\sum_{j=1}^Ma_j(R)u_j^-));L_2(\Gamma^R)\|=
O(e^{-\gamma R}),\quad R\rightarrow\infty.
\end{equation}
In view of
$$
(\mathcal{N}+i\zeta\mathcal{D})Z_l^R|_{\Gamma^R}=
(\mathcal{N}+i\zeta\mathcal{D})(u^+_l+\sum_{j=1}^Ma(R)_ju_j^-)|_{\Gamma^R},
$$
from (\ref{5.7}) it follows
\begin{equation}\label{5.8}
    \|\mathcal{N}(Z_l^R-(u^+_l+\sum_{j=1}^Ma_j(R)u_j^-));L_2(\Gamma^R)\|=
O(e^{-\gamma R}),\quad R\rightarrow\infty.
\end{equation}
Making use of  (\ref{5.7}) and (\ref{5.8}), we rewrite  (\ref{5.6})
in the form
$$
(\mathcal{N}\varphi_l,\mathcal{D}\varphi_l)_{\Gamma^R}-
    (\mathcal{D}\varphi_l,\mathcal{N}\varphi_l)_{\Gamma^R}=O(e^{-\gamma
    R}),
$$
where $\varphi_l=u^+_l+\sum a_j(R)u_j^-$. According to Proposition
\ref{p4.1}, the left-hand side is equal to
$-i(1-\sum|a_j(R)|^2)+o(1)$. Thus,
$$
\sum_{j=1}^M|a_j(R)|^2=1+o(1).\;\Box
$$

\noindent{\bf Proof of Theorem \ref{PrinTh}.} Let  $Y_l$, $Z_l^R$,
and $(a_1(R),\ldots,a_M(R))$ be the same as in Proposition
\ref{p4.4}. Substitute  $u=v=U_l:=Y_l-Z_l^R$ in the Green formula
(\ref{Green formula in GR}). Since $U_l$ satisfies the two first
equations in (\ref{problem in GR for X}), we have
\begin{equation}\label{5.9}
    (\mathcal{N}U_l,\mathcal{D}U_l)_{\Gamma^R}-
    (\mathcal{D}U_l,\mathcal{N}U_l)_{\Gamma^R}=0.
\end{equation}
We set
\begin{equation}\label{5.10}
    \varphi_l=u^+_l+\sum_{j=1}^M a_j(R)u_j^-,\quad \psi_l=u^+_l+\sum_{j=1}^M S_{lj}u_j^-
\end{equation}
and rewrite $U_l$ in the form
$$
U_l-Y_l-Z_l^R=(Y_l-\psi_l)+(\psi_l-\varphi_l)+(\varphi_l-Z_l^R).
$$
Note that  $(Y_l-\psi_l)|_{\Gamma^R}=O(e^{-\gamma R})$ by virtue of
(\ref{Y}). In view of (\ref{5.7}), (\ref{5.8}), and Proposition
\ref{p4.4}, this enables us to pass from (\ref{5.9}) to
\begin{equation}\label{5.11}
    (\mathcal{N}(\psi_l-\varphi_l),\mathcal{D}(\psi_l-\varphi_l))_{\Gamma^R}-
    (\mathcal{D}(\psi_l-\varphi_l),\mathcal{N}(\psi_l-\varphi_l))_{\Gamma^R}=O(e^{-\Lambda
    R}),
\end{equation}
where $\Lambda=\min\{\beta,\gamma\}$. The left-hand side can be
immediately calculated and is equal to

$$i\sum_{j=1}^M|a_j(R)-S_{lj}|^2+O(e^{-\beta R});$$
to see that, it suffices to use (\ref{5.10}) and Proposition
\ref{p4.1}. Finally, we obtain
$$
\sum_{j=1}^M|a_j(R)-S_{lj}|^2=O(e^{-\Lambda R}).\;\Box
$$

\end{document}